\documentclass[review]{elsarticle}
\UseRawInputEncoding 
\usepackage{amsmath}
\usepackage{amssymb}
\usepackage{amsthm}
\usepackage{psfrag} 
\usepackage{amsfonts,amsmath,enumerate,graphicx,enumerate}
\usepackage{lineno}
\modulolinenumbers[5]

\journal{Journal of \LaTeX\ Templates}

\allowdisplaybreaks

\newcommand{\Lower}[1]{\smash{\lower 1.3ex  \hbox{#1}}} 

\begin{document}

\begin{frontmatter}

\title{Using cluster analysis on municipal statistical data to configure public policies about Water, Sanitation and Hygiene in Venezuela}



\author[mymainaddress]{Sa\'ul E. Buitrago-Boret\corref{mycorrespondingauthor}}
\cortext[mycorrespondingauthor]{Corresponding author}
\ead{sbutrago@usb.ve}

\author[mymainaddress]{Roger Mart\'{\i}nez-Rivas}
\ead{rogermartinez130962@gmail.com}

\author[mymainaddress]{Josefina Fl\'orez-Diaz}
\ead{jflorez@usb.ve}

\author[mymainaddress]{Rodrigo Mijares-Seminario}
\ead{rmijares@usb.ve}

\author[mymainaddress]{Elena Rinc\'on}
\ead{erincon58@gmail.com}

\address[mymainaddress]{Universidad Sim\'on Bol\'{\i}var, Caracas, Venezuela}

\begin{abstract}
Objective: The aim of this research is to demonstrate how the use of hierarchical 
cluster analysis on 366 municipalities and other minor entities (parishes) of Venezuela, 
could be useful to consider regional differences and similarities between territorial 
entities when designing national public policies of Water, Sanitation and Hygiene 
(WASH) based on evidence. \\
Methods and results: Consider data from various sources to characterize the population 
of Venezuela through their territorial entities. Select variables at the level of the 
territorial entities to cover demographic characteristics, mortality and nutrition, 
coverage of reliable water and sanitation services, access to education, and access 
to information and communication technologies. Classify the territorial entities 
into a limited number of mutually exclusive groups using hierarchical clustering 
techniques and based on proximity in the multi-dimensional space. Adjust of assignments, 
reallocating some entities into a different group based on the specialists' opinion 
about its hierarchy in the cities regional system and its geographic location. 
Define an indicator to verify the consistency of the groups built. Conduct a 
statistical analysis to confirm separation of the groups. Demonstrate the utility 
of the results with some examples of common analysis when building a sanitary public 
policy, using seven distinct groups of recommendations depending on each cluster. \\
Conclusions: Cluster analysis can be a useful method to analyse relevant differences 
between territorial entities when designing national public policies based on evidence.
\end{abstract}


\end{frontmatter}

\nolinenumbers

\section{Introduction} \label{Sec_Introduction}

Access to adequate water and sanitation facilities does not necessarily lead 
to better health. Besides this affirmation, identify how relevant are several 
social determinants to explain differences in many water-related diseases, is 
a complex problem difficult to solve, much more when designing a public policy. 
Reducing diseases and mortality rates on waterborne and other water-related 
diseases is due not only to inadequate quality and quantity of water for human 
consumption, but also to poor hygiene practices that are commonly associated 
with poverty, bad access to education and adequate information, among other 
determinants.

Even within the same nation, there are differences in social determinants by 
region, and between urban and rural areas (see~\cite{R1,R2,R3}). Consequently, it is 
essential to perform analyses that discover regional differences, to design 
strategies that promote a reduction of water-borne diseases adapted to each 
context.

The political, economic, and social crisis that Venezuela is going through, 
and the weakening of the competent public institutions have negatively affected 
WASH services, both at the level of infrastructure and service provision, to 
which are added inadequate hygienic practices at home (see~\cite{R4,R5,R6,R7,R8,R9}).

Based on the most recent statistical data available at the level of federal 
entities and municipalities, an epidemiological profile of Venezuela was prepared 
regarding some water-borne diseases (WBD), demographic and socioeconomic indicators 
were described, the coverage of water services and sanitation, the assistance 
to educational services, and the access to information and communication technologies 
was analysed. These analyses made it possible to indicate the situation of the 
country and reflect regional differences, which were categorized through a 
hierarchical cluster analysis.

The politic-administrative division of Venezuela covers 23 states, 1 Capital District 
and 1 federal dependencies. Internally each state is divided into municipalities. 
The country has 335 municipalities that are also divided in parishes.

Data / information processing and analysis from various sources was carried out 
with the aim of characterizing the population of Venezuela. To process the data, 
at the states level of Venezuela, a database was created, and some basic statistical 
analyses were carried out, which led to tables of total data; this to graphically 
visualize some characteristics of the territorial entities.

\begin{figure}[htp!]
   \begin{center}
   \includegraphics[width=7.5cm]{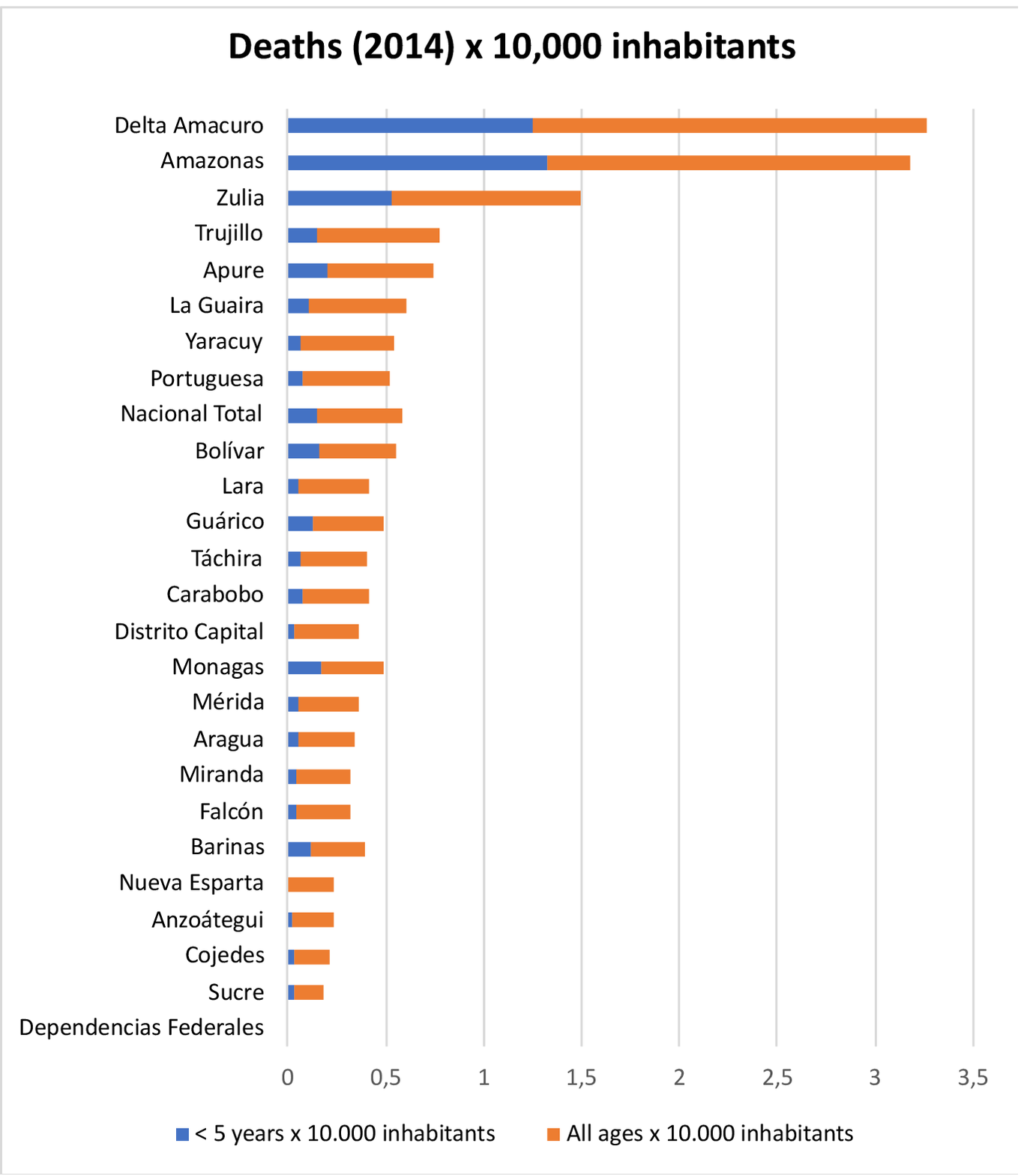}
   \includegraphics[width=7.5cm]{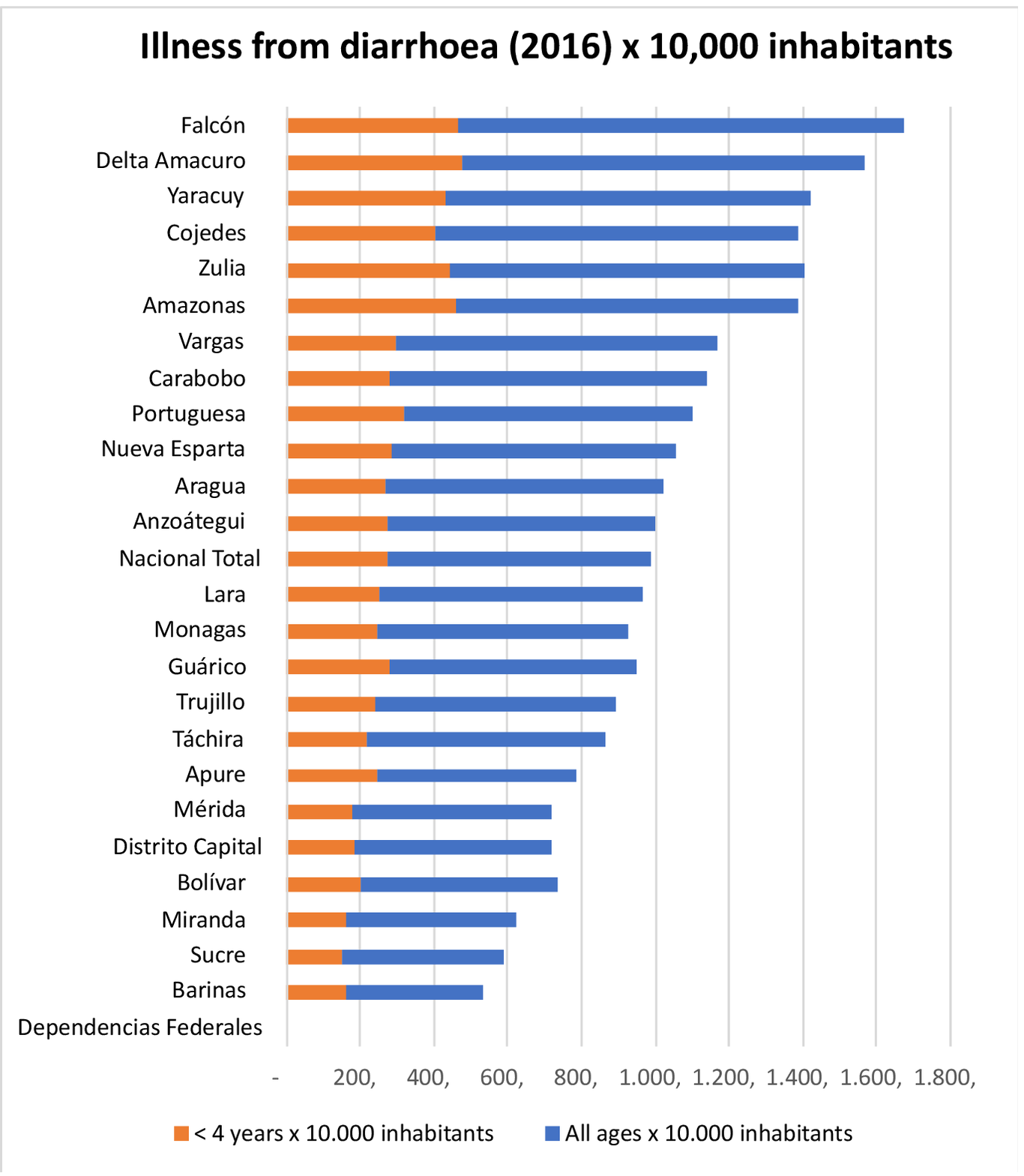}
   \vspace*{-2mm}
   \caption {On the top, deaths from diarrhoea and gastroenteritis of presumed 
    infectious origin x 10.000 inhabitants, by federal entity. Venezuela. 2014. 
    See~\cite{R10,R11}. On the bottom, illness from diarrhoea in children less than 4 
    years of age and all ages x 10.000 inhabitants, by federal entity. Venezuela. 
    2016. See [10,12].}
   \label{deaths-from-diarrhoea-and-gastroenteritis}
   \end{center}
\end{figure}

Fig.~\ref{deaths-from-diarrhoea-and-gastroenteritis} (top) shows data related to 
deaths from diarrhoea and gastroenteritis of presumed infectious origin registered 
in Venezuela in 2014 (see~\cite{R10,R11}); it illustrates the relative importance 
of deaths with respect to the total population of each federal entity.

The higher prevalence of Delta Amacuro, Amazonas, and Zulia with respect to the 
national total can be observed in the registry of deaths associated with the ingestion 
of water and food in unsafe conditions. In these three states, the prevalence is 
double the national average. It can also be noted that the participation of children 
under 5 years of age in the composition of the general rate of deaths due to diarrhoea 
and gastroenteritis in these three entities accounted for more than half of the cases, 
which tells us about the precariousness of water, sanitation and hygiene, limitations 
in vaccination, limitation in child deworming campaigns, precariousness of the shelter 
and the socioeconomic situation of the families, among other possible reasons.

The states Trujillo, Apure, Vargas, Yaracuy, and Portuguesa, they also register 
prevalence above the national average. In these cases, although the population group 
of children under 5 years of age has a lower incidence, the conditions of access to 
water and hygiene in homes must also be weak.

Fig.~\ref{deaths-from-diarrhoea-and-gastroenteritis} (bottom) shows 2016´s records of 
illness from diarrhoea in children less than 4 years of age and all ages by 10 thousand 
of inhabitants (see~\cite{R10,R12}). According to this data, the number of sicks is higher 
than deaths, remaining a high incidence of children under 4 years, and a major prevalence 
of Delta Amacuro, Amazonas, and Zulia adding Falc\'on, Yaracuy and Cojedes over the rest 
of federal entities in Venezuela.

In this research, the objects under study will be the 335 municipalities that make up 
Venezuela, excluding the Libertador Municipality (the only municipality of the Distrito 
Capital) for which the 22 parishes that comprise compound it will be considered, 
in the same way, we will proceed with the Vargas municipality (the only municipality in 
the Vargas state), considering its 11 parishes. Thus, ending with 366 objects for the study. 
According to the Living Conditions Survey ENCOVI 2020 (see~\cite{R13,R14,R15,R16,R17,R18}), 
the country registered in 2020 a declining population of the order of 28.5 million 
inhabitants, distributed in 366 municipalities and parishes referred above. 
The municipalities and the parishes present notable differences regarding population size, 
percentage of urban and rural population, levels of human development, coverage of water 
and sanitation services, presence of indigenous communities and informal settlements, 
among other key attributes to identify the most convenient type of intervention.

The idea pursued in this research is to demonstrate how the use of hierarchical cluster 
analysis on territorial entities of Venezuela, could be useful to consider regional 
differences and similarities between territorial entities when designing national 
public policies of Water, Sanitation and Hygiene (WASH) based on evidence.

\section{Methods} \label{Sec_Methods}

To process the data from the ENCOVI 2020 survey, at the municipalities of Venezuela 
and parishes of Distrito Capital and Vargas state, a database was created.

Regarding data on the coverage of WASH and other services at homes, the ENCOVI survey 
(see~\cite{R13,R14,R15,R16,R17,R18}) includes 26 attributes / variables at the level 
of territorial entities; of these, 15 were selected ($x_1$ to $x_{15}$) within the 
following 6 areas:
\begin{enumerate}
  \item coverage of reliable water and sanitation services
  \begin{itemize}
    \item Percentage of households with safely managed drinking water supplies ($x_1$).
    \item Percentage of households using safely managed sanitation services ($x_2$).
  \end{itemize}
  \item health status of the population
  \begin{itemize}
    \item Mortality rate for children under 5 years of age ($x_3$).
    \item Percentage of children under 5 years of age who are underweight ($x_4$).
  \end{itemize}
  \item socioeconomic data
  \begin{itemize}
    \item Percentage of households living below the extreme poverty line ($x_5$).
    \item Total employment rate ($x_6$).
  \end{itemize}
  \item demographic situation
  \begin{itemize}
    \item Percentage of population ($x_7$).
    \item Percentage of adolescent mothers ($x_8$).
    \item Demographic dependency ratio ($x_9$).
  \end{itemize}
  \item access to education in public institutions
  \begin{itemize}
    \item Rate of attendance at a teaching centre for children from 3 to 5 years old ($x_{10}$).
    \item Rate of attendance at a teaching centre for children from 6 to 11 years old ($x_{11}$).
    \item Rate of attendance at a teaching centre for adolescents aged 12 to 17 years ($x_{12}$).
    \item Percentage of the population aged 15 to 64 that reaches at least complete primary education ($x_{13}$).
  \end{itemize}
  \item access to information and communication technologies
  \begin{itemize}
    \item Percentage of households that have computers ($x_{14}$).
    \item Percentage of households that have internet access ($x_{15}$).
  \end{itemize}
\end{enumerate}

Next, new variables, $z_1$ to $z_{15}$, were defined as follows. The variation ranges 
of each attribute $x_i$ were reduced to the intervals defined by their minima and maxima, 
and then normalized to the range 0 to 100, for the population of the 366 territorial 
entities, as shown in equation
\footnote{It is also known as the min-max normalization method ``Min-Max: Min-Max normalization 
is the process of taking data measured in its engineering units and transforming it to a 
value between 0 and 1. Whereby the lowest (min) value is set to 0 and the highest (max) 
value is set to 1. This provides an easy way to compare values that are measured using different 
scales or different units of measure. The normalized value is defined as: 
$MM(x_{ij})=(x_{ij}-x_{min})/(x_{max}-x_{min})$'', (see~\cite{R19}).} \ref{E1}.

\begin{equation}
\label{E1}
z_i = \frac{x_i - \min_{1 \leq j \leq 15} x_j}{\max_{1 \leq j \leq 15} x_j - \min_{1 \leq j \leq 15} x_j} 100, 
\quad {\rm for} \ 1 \leq i \leq 15.
\end{equation}

Some basic statistical analyses were carried out, which led to tables of total data 
and radial charts; this to graphically visualize some characteristics of the territorial 
entities.

When analysing the data, it is evident that within the territorial entities there are 
relevant differences regarding the six areas named before. 

To develop a water, sanitation and hygiene strategy adapted to these regional differences, 
a hierarchical cluster analysis was carried out for the 366 territorial entities to group 
them into relatively homogeneous categories. The idea is to discover associations and structures 
in the data that are not obvious a priori, but that can be useful once they have been found. 
The objective of this type of analysis is not to statistically explain the phenomenon, but 
to classify some objects (in this case territorial entities defined by common characteristics), 
into a limited number of mutually exclusive groups.

The goal of cluster analysis is to group objects into clusters such that objects within the 
same cluster are similar in some sense, and objects from different clusters are dissimilar 
(see~\cite{R20,R21}). The hierarchical clustering algorithms selected proceed agglomerative 
(bottom-up). In agglomerative algorithms, each object starts as a singleton cluster and the 
``clusters'' are merged successively. 
Most clustering methods form clusters based on the proximity between data points in the 
multi-dimensional space. A commonly used measure of proximity between a pair of data points 
(objects) is the squared Euclidean Distance. Hierarchical algorithms output a rooted tree 
structure that can be represented as a dendrogram. All analyses are conducted using Statistics 
Toolbox of MATLAB (see~\cite{R22}).

\begin{figure}[htp!]
   \begin{center}
   \includegraphics[width=12cm]{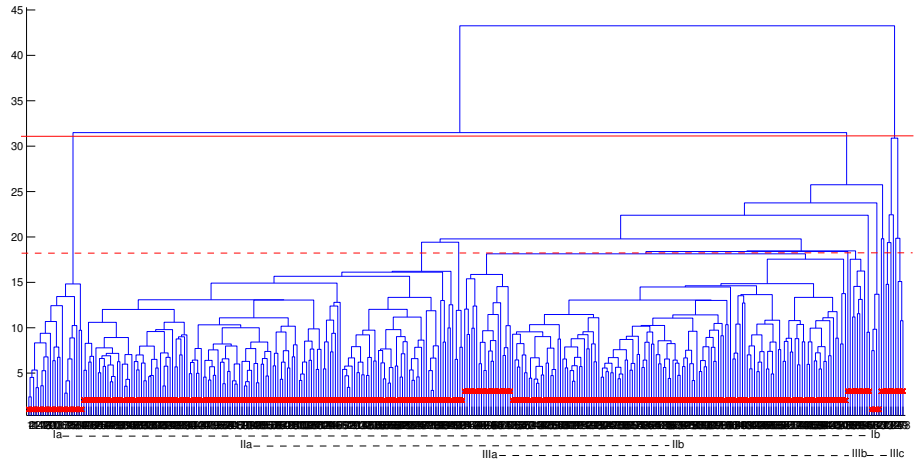}
   \vspace*{-2mm}
   \caption {Dendrogram associated with the cluster analysis for the 366 territorial entities. Year 2020.}
   \label{Dendrogram-366-entities}
   \end{center}
\end{figure}

The hierarchical dendrogram corresponding to the analysis carried out, using a solid red cut 
line, is shown in Fig.~\ref{Dendrogram-366-entities}, where we can see 3 large very well 
differentiated groups. The group on the left with 23 territorial entities of the Distrito 
Capital and the States of Anzoátegui, Carabobo, M\'erida, Miranda, Nueva Esparta, and T\'achira. 
The central group that contains most of the country's territorial entities (335 between 
municipalities and parishes). And the group on the right with 8 municipalities in the states 
of Amazonas, Barinas, and Delta Amacuro.

With this general overview, it was decided to review in detail the composition of each group. 
Carrying out a group cut further down (shown by the dash red line on the dendrogram in 
Fig.~\ref{Dendrogram-366-entities} produces 16 subgroups, where the original group on the 
left remains the same, the original in the centre is divided into 10 subgroups and the 
original on the right is divided into 5 subgroups.

Next, based on these 16 new subgroups, the municipalities were classified into three large 
categories, relocating some subgroups of municipalities that, according to a taxonomic analysis, 
belonged to the original central group, but in the opinion of the specialists they would be 
better located in some of the extreme original groups, if other variables were considered, 
such as hierarchy in the regional system of cities and remote geographic location with 
respect large metropolitan areas. The illustration below the dendrogram in 
Fig.~\ref{Dendrogram-366-entities} shows the 7 subcategories Ia (cluster 16), Ib (cluster 15), 
IIa (clusters 7 and 8), IIb (cluster 2), IIb (cluster 2), IIIa (clusters 3, 4, 8, 10 and 13), 
IIIb (cluster 1), and IIIc (clusters 5, 6, 11, 12, and 14).

The characteristics of each category are described below. \\
Category I: Municipalities and parishes that represent jurisdictions within metropolitan 
areas that are generally well served with respect to the coverage of water and sanitation 
services, with favourable indicators regarding their health, demographic, and socioeconomic 
situation, with high percentages of households with access to information and communication 
technologies. \\
Category II: Territorial entities in intermediate and smaller cities, some capitals of 
federal entities (the vast majority of the country's municipalities), with less favourable 
indicators, but distant from the other categories in most of the indicators used. \\
Category III: Corresponds to municipalities located in remote areas, far from main populated 
centres; It comprises mostly border municipalities, with high infant mortality rates, 
unfavourable indicators regarding their demographic and socioeconomic situation, with great 
limitations in access to educational services and information and communication technologies.

In this sense, the result of the cluster analysis is not rigorously followed, but it does 
allow a general understanding of the differences that justify three categories for a future 
intervention strategy: two extreme categories and a large central majority group.

Next, the NL2 ($L_2$ norm) indicator is defined in equation (2) based on the Euclidean 
norm with weights, where certain attributes are given greater relevance, such as: the 
percentage of households that have safe managed drinking water supply services ($x_1$), 
the mortality rate of children under 5 years of age ($x_3$), the percentage of households 
that live below the extreme poverty line ($x_5$), the percentage of adolescent mothers ($x_8$), 
and the percentage of the population aged 15 to 64 that reaches at least complete primary 
school ($x_{13}$). The other attributes considered in the definition of NL2 is the denial 
of the meaning of the attributes: mortality rate for children under 5 years of age ($x_3$), 
percentage of children under 5 years of age who are underweight ($x_4$), percentage of 
households living below the extreme poverty line ($x_5$) and percentage of adolescent 
mothers ($x_8$), to have positive contributions to the indicator.

\begin{equation}
\label{E"}
NL2 = (\sum_{i=1}^{15} x_i^2 w_i)^{1/2}, \quad {\rm with} \ \sum_{i=1}^{15} w_i = 1.
\end{equation}

\begin{figure}[htp!]
   \begin{center}
   \includegraphics[width=12cm]{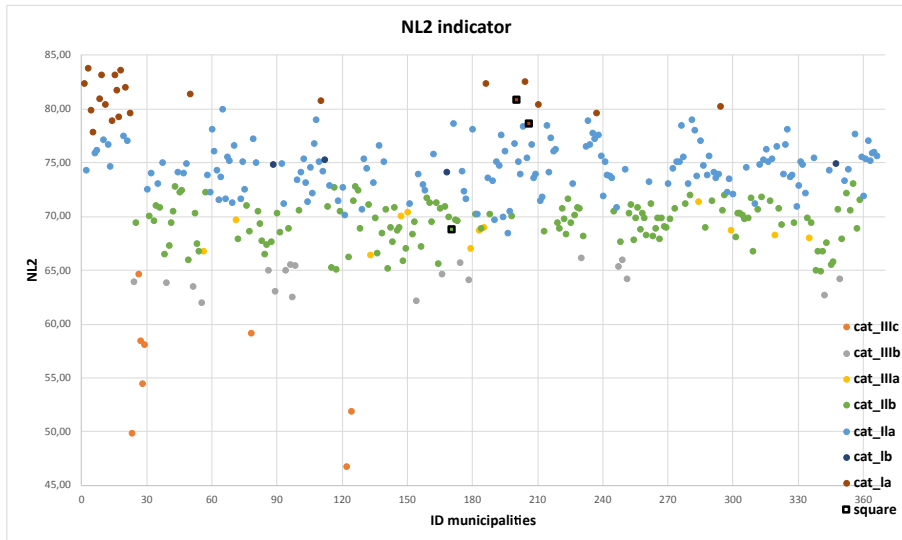}
   \vspace*{-2mm}
   \caption {Territorial entities distributed according to the indicator NL2 (L2 norm). Year 2020.}
   \label{NL2-indicator}
   \end{center}
\end{figure}

Fig.~\ref{NL2-indicator} shows how the 366 territorial entities are distributed 
according to the indicator NL2. We can appreciate how the 7 subcategories are 
separated. The abscissa axis corresponds to the number assigned to each territorial 
entity and the ordinate axis the calculated value of NL2 for each territorial entity.

A way to verify the validity and accuracy of the results is to analyse and compare 
characteristics of specific local entities. The municipalities of Moran (Lara State), 
Baruta and El Hatillo (Miranda State) have been chosen for detailed comparisons.

\begin{table}[htp!]
\caption{Values of variables for Moran, Baruta and El Hatillo Municipalities. Year 2020.}
\label{ValuesOfVariables}
\begin{center}
\resizebox{\columnwidth}{!}{
\begin{tabular}{|c|c|c|c|}
\hline
Variables        & Moran & Baruta & El Hatillo \\
                 & ID: 170 & ID: 200 & ID: 206 \\
                 & NL2: 68.81 & NL2: 80.84 & NL2: 78.66 \\
\hline
Percentage of households with safely managed drinking water supplies ($x_1$)  & 77,77\% & 91,24 \% & 79,14\% \\
Percentage of households using safely managed sanitation services ($x_2$)  & 69,80\% & 99,90 \% & 98,39\% \\
Mortality rate for children under 5 years of age ($x_3$)  & 24,85\% & 10,64\% & 10,90\% \\
Percentage of children under 5 years of age who are underweight ($x_4$)  & 9,06\% & 3\% & 3,19\% \\
Percentage of households living below the extreme poverty line ($x_5$)  & 67\% & 47\% & 51\% \\
Total employment rate ($x_6$)  & 54\% & 54\% & 50\% \\
Percentage of population ($x_7$)  & 0,44\% & 1,17\% & 0,29\% \\
Percentage of adolescent mothers ($x_8$)  & 15,05\% & 5,56\% & 6,39\% \\
Demographic dependency ratio ($x_9$)  & 59\% & 39\% & 41\% \\
Rate of attendance at a teaching centre for children from 3 to 5 years old ($x_{10}$)  & 49,37\% & 74,75\% & 74,59\% \\
Rate of attendance at a teaching centre for children from 6 to 11 years old ($x_{11}$)  & 94,76\% & 95,20\% & 92,05\% \\
Rate of attendance at a teaching centre for adolescents aged 12 to 17 years ($x_{12}$)  & 77,33\% & 90,17\% & 86,13\% \\
Percentage of the population aged 15 to 64 with at least primary education ($x_{13}$)  & 76,53\% & 95,43\% & 94,45\% \\
Percentage of households that have computers ($x_{14}$)  & 12,15\% & 50,20\% & 46,07\% \\
Percentage of households that have internet access ($x_{15}$)  & 7,28\% & 47,38\% & 42,50\% \\
Cluster (Capital number: Category; Small letter: Subcategory)  & IIb & Ia & Ia \\
\hline
\end{tabular}
}
\end{center}
\end{table}

As one can observe in Table~\ref{ValuesOfVariables}, in general the values of 
Hatillo and Baruta are relatively closer between them than when are compared with 
Moran values. However, some indicators are too close as to explain differences 
between entities (for example: variables $x_6$ and $x_{11}$), and in some cases 
values are distributed in a way that a single indicator considered alone could 
justify a different classification (for example: variable $x_7$). These simple 
observations allow to asseverate the importance of consider a wide set of variables 
and the relevance of analyse them, trying to smooth small differences and to discover 
a convincing grouping of entities.

With this idea in mind, the municipalities Mor\'an (ID:170, NL2:68.81) belonging 
to subcategory IIb, and Baruta (ID:200, NL2:80.84) and El Hatillo (ID:206, NL2:78.66) 
belonging to category Ia, are highlighted with a black square in Fig.~\ref{NL2-indicator}

It can be noted that municipalities Baruta and El Hatillo in Miranda state, part of 
the metropolitan area of Caracas, are clearly far from Mor\'an in Lara state and from 
other entities categorized in II (compare red and green dots in Fig.~\ref{NL2-indicator}). 
Indeed, according with further analysis, these two first entities have a more favourable 
situation than Mor\'an, regarding services available and hygiene practices at home.

Furthermore, it is convenient to advice that within each municipality and parish there 
may also appear some differences, especially in urban areas, that can justify detailed 
analysis for future interventions. For example, in the El Hatillo municipality of Miranda 
state, a less favourable situation can be observed in the rural settlement of Gavil\'an, 
an inner community, in comparison with El Hatillo town, its capital. The same situation 
appears with the neighbourhoods of La Palomera and Ojo de Agua, slum settlements inside 
Municipality of Baruta, compared to residential areas of higher income groups. In general, 
the communities located in the urban peripheries have situations that are more like rural 
areas or small populated rural centres of the country.

\section{Results} \label{Sec_Results}

The municipalities and parishes of category type I comprise 27 jurisdictions, presented 
by 14 parishes of the Distrito Capital and 13 municipalities located in the states of 
Anzo\'ategui, Bol\'{\i}var, Carabobo, Lara, M\'erida, Miranda, Nueva Esparta, T\'achira, 
and Zulia. The estimated population in 2020 exceeds 6 million 400 thousand inhabitants 
(22,46\% of the total population of the country). The health conditions, water and 
sanitation coverage, economic, demographic, educational and services situation in these 
jurisdictions are the most favourable in the national context. They are located within 
large metropolitan areas, which ensures access to services, although within them there 
may be informal settlements with unfavourable situations, as in the cases of Baruta and 
El Hatillo in the State of Miranda. The data analysis distinguishes two subcategories, 
Ia (14 parishes and 9 municipalities) and Ib (4 municipalities).

The municipalities and parishes of category type II comprises 299 jurisdictions, 
consisting of eight parishes of the Distrito Capital and 291 municipalities located in 
all federal entities in the country. The estimated population in 2020 reaches about 21 
million inhabitants (73,68\% of the total population of the country). Their health 
conditions, water and sanitation coverage, socioeconomic, demographic, educational and 
services situation is around the national average. The parishes of the Distrito Capital 
and several municipalities are in large metropolitan areas, but the majority belong to 
intermediate cities and small rural populations. The data analysis distinguishes two 
subcategories, IIa (159 entities) and IIb (140 entities).

The municipalities and parishes of category type III comprise 40 municipalities located 
in Amazonas, Anzoátegui, Apure, Aragua, Barinas, Bol\'{\i}var, Delta Amacuro, Falc\'on, 
Gu\'arico, Lara, M\'erida, Monagas, Portuguesa, T\'achira, Trujillo, and Zulia. The 
estimated population in 2020 reaches around 1 million 100 thousand inhabitants (3,86\% of 
the total population of the country). Their health conditions, water and sanitation 
coverage, socioeconomic, demographic, educational and services situation are worse than 
the national average (depending on the attribute considered). They correspond to small 
rural populations distant from the main cities of the regions to which they belong, so 
their access to goods and services is limited. The data analysis distinguishes three 
subcategories, IIIa (12 municipalities), IIIb (20 municipalities) and IIIc (8 municipalities).

\begin{figure}[htp!]
   \begin{center}
   \includegraphics[width=12cm]{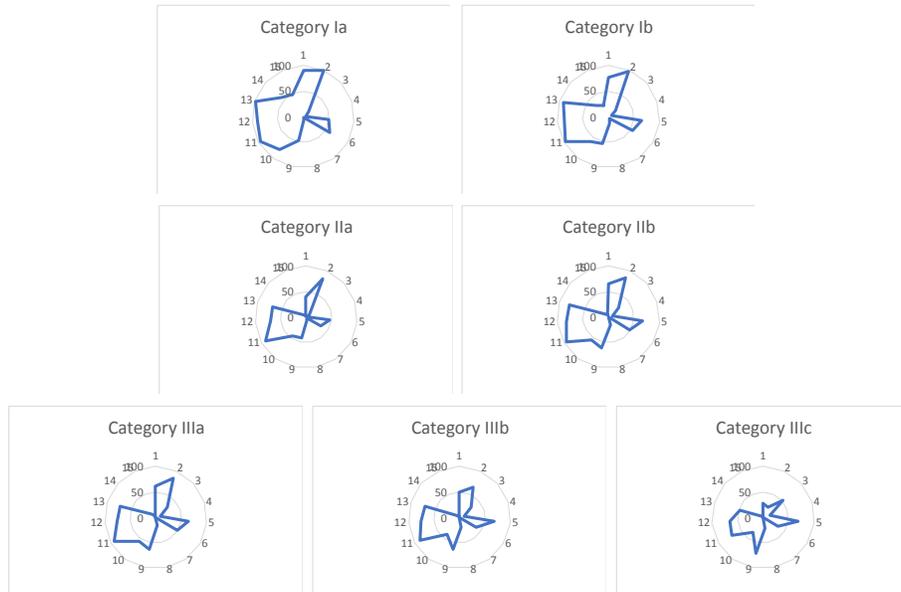}
   \vspace*{-2mm}
   \caption {Radial charts of average values in the jurisdictions for each of the seven subcategories. Year 2020.}
   \label{RadialCharts}
   \end{center}
\end{figure}

Table~\ref{AverageValuesAttributes} shows the average values for the attributes, and Fig.~\ref{RadialCharts}
shows the radial charts for each of the subcategory belonging to categories I, II and III.

\begin{table}[htp!]
\caption{Average values of the 15 attributes in the jurisdictions for each of the seven subcategories. Year 2020.}
\label{AverageValuesAttributes}
\begin{center}
\resizebox{\columnwidth}{!}{
\begin{tabular}{|c|c|c|c|c|c|c|c|}
\hline
Attributes / Categories & Ia & Ib & IIa & IIb & IIIa & IIIb & IIIc \\
\hline
Percentage of households with safely managed drinking water supplies ($x_1$)  & 90,61 & 76,80 & 40,49 & 66,41 & 61,79 & 49,02 & 27,42 \\
Percentage of households using safely managed sanitation services ($x_2$)  & 99,60 & 97,48 & 81,78 & 84,52 & 85,63 & 65,40 & 22,79 \\
Mortality rate for children under 5 years of age ($x_3$)  & 14,98 & 18,63 & 10,74 & 27,60 & 31,65 & 29,66 & 51,43 \\
Percentage of children under 5 years of age who are underweight ($x_4$)  & 4,09 & 6,45 & 4,22 & 9,56 & 10,01 & 10,12 & 14,03 \\
Percentage of households living below the extreme poverty line ($x_5$)  & 49,36 & 64,50 & 47,00 & 67,66 & 64,33 & 68,80 & 69,63 \\
Total employment rate ($x_6$)  & 59,43 & 54,00 & 55,52 & 47,84 & 64,33 & 38,75 & 34,50 \\
Percentage of population ($x_7$)  & 0,33 & 3,76 & 0,33 & 0,15 & 0,07 & 0,14 & 0,04 \\
Percentage of adolescent mothers ($x_8$)  & 7,03 & 12,78 & 7,32 & 16,20 & 14,78 & 19,94 & 21,51 \\
Demographic dependency ratio ($x_9$)  & 46,43 & 54,25 & 41,00 & 61,14 & 62,42 & 63,75 & 71,75 \\
Rate of attendance at a teaching centre for children from 3 to 5 years old ($x_{10}$)  & 78,15 & 60,62 & 44,60 & 55,98 & 56,39 & 40,64 & 35,02 \\
Rate of attendance at a teaching centre for children from 6 to 11 years old ($x_{11}$)  & 37,07 & 96,64 & 90,90 & 95,80 & 94,01 & 89,94 & 70,15 \\
Rate of attendance at a teaching centre for adolescents aged 12 to 17 years ($x_{12}$)  & 89,98 & 85,73 & 70,50 & 81,89 & 74,89 & 76,10 & 65,75 \\
Percentage of the population aged 15 to 64 with at least primary education ($x_{13}$)  & 97,14 & 92,85 & 69,04 & 80,29 & 72,90 & 70,99 & 47,42 \\
Percentage of households that have computers ($x_{14}$)  & 57,80 & 32,50 & 7,12 & 11,24 & 9,99 & 8.08 & 3,34 \\
Percentage of households that have internet access ($x_{15}$)  & 48,63 & 23,40 & 3,43 & 5,58 & 3,37 & 3,46 & 0,81 \\
\hline
\end{tabular}
}
\end{center}
\end{table}

\begin{figure}[htp!]
   \begin{center}
   \includegraphics[width=12cm]{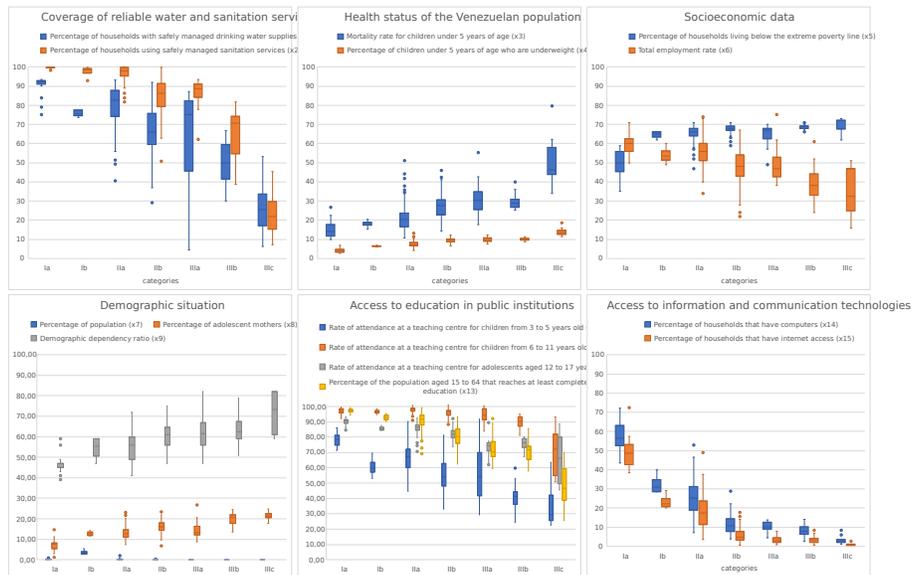}
   \vspace*{-2mm}
   \caption {Box plot regarding data on the coverage of WASH and other services at homes 
         in the areas and attributes for the seven subcategories. Year 2020.}
   \label{BoxPlots}
   \end{center}
\end{figure}

Fig.~\ref{BoxPlots} presents the box plot
\footnote{WASH and other services at homes outcomes are presented as box plots, 
which represent the middle 50\% of the data, ranging from the upper boundary 
(75th percentile) to the lower boundary (25th percentile). Box lines indicate 
median values. Vertical lines extending from the box indicate minimum and maximum 
values and dots are outliers.} regarding data on the coverage of WASH and 
other services at homes in the six areas defined before (see section Materials and 
methods) for the 15 attributes for the seven subcategories coming from the cluster 
analysis performed. It can be observed how the seven subcategories separate to 
each other for each of the 15 attributes.

\begin{figure}[htp!]
   \begin{center}
   \includegraphics[width=12cm]{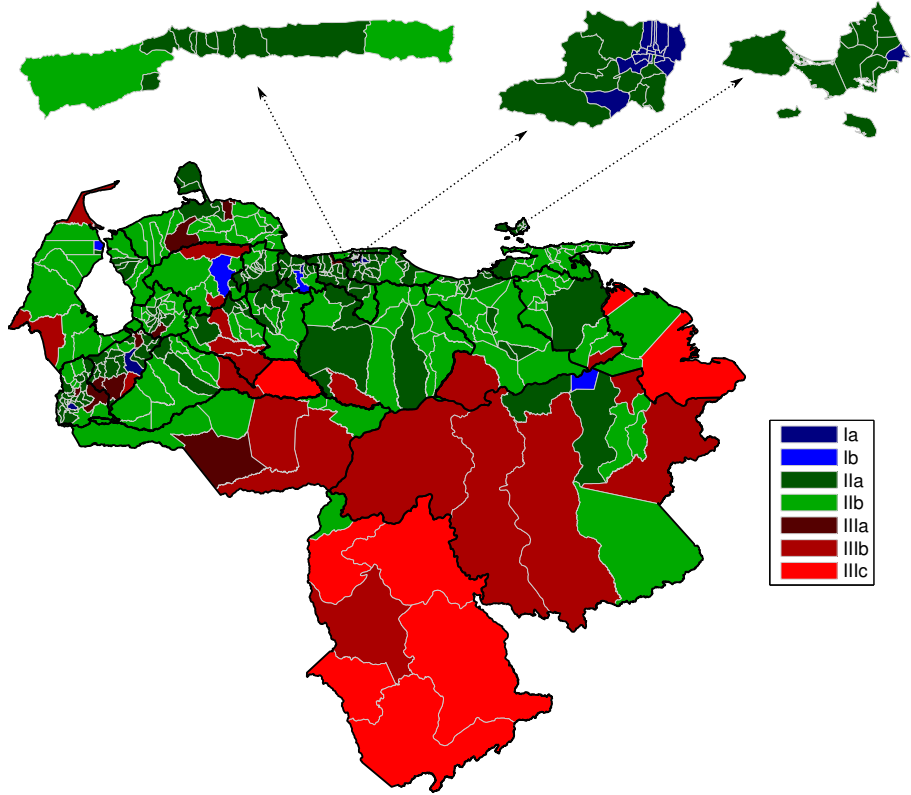}
   \vspace*{-2mm}
   \caption {Municipalities and parishes of Venezuela under category Ia, Ib, IIa, IIb, IIIa, 
      IIIb and IIIc. Top left, Vargas state with 11 parishes. Top centre, Distrito Capital 
      with 22 parishes. Top right, Nueva Esparta state with 11 municipalities. Bottom, 
      Venezuela with 335 municipalities excluding Distrito Capital and the Vargas state.}
   \label{mapa-venezuela}
   \end{center}
\end{figure}

Fig.~\ref{mapa-venezuela} shows the 366 territorial entities distributed 
according to the 7 subcategories obtained by the cluster and taxonomic analysis.

As an example of how a clustering analysis allows design strategies to perform a water and 
sanitation program of intervention, Table~\ref{DesignStrategies} below shows how different 
types of intervention can be more or less suitable to different type of municipalities, 
according to their attributes. The taxonomy of municipalities allows a quick analysis of 
viability once the universe of local entities has been well categorized.

\begin{table}[htp!]
\caption{Example of different types of intervention to different type of municipalities, 
according to their attributes. Year 2020. Legend: 5 more suitable / reliable, 4 relatively 
suitable / reliable, 3 less reliable / less suitable, 2 unuseful / unapplicable, 
1 need more information.}
\label{DesignStrategies}
\begin{center}
\resizebox{\columnwidth}{!}{
\begin{tabular}{|c|c|c|c|c|c|c|c|}
\hline
Examples of public policies on Water, Sanitation and Hygiene / Categories & Ia & Ib & IIa & IIb & IIIa & IIIb & IIIc \\
\hline
Designing / building / improving large and medium water \& sanitation infrastructure systems  & 5 & 5 & 5 & 5 & 3 & 2 & 2 \\
Designing / building / improving small and domestic water \& sanitation systems  & 3 & 3 & 3 & 3 & 5 & 5 & 5 \\
Providing information to ensure safe water, sanitation, and hygiene practices at household level through workshops & 4 & 4 & 4 & 4 & 5 & 5 & 5 \\
Providing information to ensure safe water, sanitation and hygiene practices at household level using TICs & 5 & 5 & 4 & 3 & 2 & 2 & 2 \\
Humanitarian assistance and donations on WASH services & 1 & 1 & 1 & 1 & 5 & 5 & 5 \\
\hline
\end{tabular}
}
\end{center}
\end{table}

\section{Conclusions} \label{Sec_Conclusions}

Venezuela has a great diversity of urban and rural contexts, with very unequal 
infrastructure and WASH services conditions, additionally, there are great cultural 
and geographical differences that must be considered in the application of strategies 
for the promotion of best WASH practices at household level.

Regarding our original objective, to characterize the population of Venezuela, 
we find strong support for the idea that municipalities clusters in Venezuela 
are associated with households with safely managed drinking water supplies, and 
using safely managed sanitation services, mortality rate for children under 5 
years of age, children under 5 years of age who are underweight, households living 
below the extreme poverty line, employment rate, population rate, adolescent mothers, 
demographic dependency, attendance at a teaching centre for children, and for 
adolescents, population aged 15 to 64 that reaches at least complete primary, 
households with computers, and internet access.

Cluster analysis of the most recent statistical data available at the level of 
municipalities of Venezuela related to WASH and other services at homes can segment 
the municipalities into seven distinct groups, Ia, Ib, IIa, IIb, IIIa, IIIb, and IIIc, 
for unique recommendations in the future in the field of Water, Sanitation and Hygiene 
interventions.

The entities belonging to category I correspond to better served municipalities and 
parishes, with more benevolent indicators, belonging to main metropolitan areas. 
The entities in category II are municipalities in intermediate and smaller cities, 
some capitals of federal entities, with less favourable indicators. Finally, the 
entities in category III are remote municipalities, far from metropolitan areas, 
with unfavourable indicators in almost all aspects.

This segmentation of groups is relevant when it is planning a humanitarian WASH 
Intervention at national and subnational levels. Close municipalities with a menu 
of options to provide information to the people, with several community facilities 
near to vulnerable communities are more affordable to plan the distribution of human, 
material and technical resources than those far locations where exist big gaps on water, 
sanitation and hygiene conditions, with historical disadvantages on socioeconomic and 
demographic indicators, and where transportation, points of information and availability 
of resources are hard to overcome.

Despite some potential limitations, this approach yields valuable intelligence to inform 
local government planning service and at the same time offers great potential for further 
research into its use in informing possible preventative WASH services.

\section*{Abbreviations}
\noindent
BM: Banco Mundial; ENCOVI: Encuesta Nacional de Condiciones de Vida; 
INE: Instituto Nacional de Estadística; OCHA: UN Office for the Coordination of 
Humanitarian Affairs; PAHO: Pan American Health Organization; WASH: Water, Sanitation 
and Hygiene; WHO: World Health Organization.



\bibliography{cluster-BuitragoMartinezFlorezMijaresRincon}

\end{document}